\documentclass[12pt,reqno]{amsart}
\usepackage{amssymb,enumerate,ifthen}
\allowdisplaybreaks[1]
\def\N{\mathbb{N}}
\def\R{\mathbb{R}}
\def\Cdot{\!\cdot\!}

\newtheorem{theorem}{Theorem}
\newtheorem*{theorem*}{Theorem}
\def\Thm#1#2{\ifthenelse{\equal{#1}{*}}{\begin{theorem*}#2\end{theorem*}}
             {\begin{theorem}\label{T#1}#2\end{theorem}}}
\newtheorem{Atheorem}{Theorem}

\def\thm#1{Theorem~\ref{T#1}}

\newtheorem{corollary}[theorem]{Corollary}
\newtheorem*{corollary*}{Corollary}
\def\Cor#1#2{\ifthenelse{\equal{#1}{*}}{\begin{corollary*}#2\end{corollary*}}
             {\begin{corollary}\label{C#1}#2\end{corollary}}}

\newtheorem{lemma}[theorem]{Lemma}
\newtheorem*{lemma*}{Lemma}
\def\Lem#1#2{\ifthenelse{\equal{#1}{*}}{\begin{lemma*}#2\end{lemma*}}
             {\begin{lemma}\label{L#1}#2\end{lemma}}}
\def\lem#1{Lemma~\ref{L#1}}
\newtheorem{Alemma}{Lemma}

\def\eq#1{{\rm(\ref{E#1})}}
\def\Eq#1#2{\ifthenelse{\equal{#1}{*}}
  {\begin{equation*}\begin{aligned}#2\end{aligned}\end{equation*}}
  {\begin{equation}\begin{aligned}\label{E#1}#2\end{aligned}\end{equation}}}

\def\diag{\mathop{\hbox{\rm diag}}\nolimits}

\begin{document}
\vspace{5mm}

\date{\today}

\title[The invariance of arithmetic mean to generalized Bajraktarevi\'c means]{On the invariance of the arithmetic mean with respect to generalized Bajraktarevi\'c means}

\author[R.\ Gr\"unwald]{Rich\'ard Gr\"unwald}
\address[R.\ Gr\"unwald]{Doctoral School of Mathematical and Computational Sciences, University of Debrecen, H-4002 Debrecen, Pf. 400, Hungary}
\email{richard.grunwald@science.unideb.hu}

\author[Zs. P\'ales]{Zsolt P\'ales}
\address[Zs. P\'ales]{Institute of Mathematics, University of Debrecen, H-4002 Debrecen, Pf. 400, Hungary}
\email{pales@science.unideb.hu}

\thanks{The research of the first author was supported by the \'UNKP-20-3 New National Excellence Program of the Ministry of Human Capacities.\ The research of the second author was supported by the K-134191 NKFIH Grant and the 2019-2.1.11-TÉT-2019-00049 and the EFOP-3.6.1-16-2016-00022 projects. The last project is co-financed by the European Union and the European Social Fund.}

\begin{abstract}
The purpose of this paper is to investigate the following invariance equation involving two $2$-variable generalized Bajraktarevi\'c means, i.e., we aim to solve the functional equation
$$
	f^{-1}\bigg(\frac{p_1(x)f(x)+p_2(y)f(y)}{p_1(x)+p_2(y)}\bigg)+g^{-1}\bigg(\frac{q_1(x)g(x)+q_2(y)g(y)}{q_1(x)+q_2(y)}\bigg)=x+y \qquad(x,y\in I),
$$
where $I$ is a nonempty open real interval and $f,g:I\to\mathbb{R}$ are continuous, strictly monotone and $p_1,p_2,q_1,q_2:I\to\mathbb{R}_+$ are unknown functions.
The main result of the paper shows that, assuming four times continuous differentiability of $f$, $g$, twice continuous differentiability of $p_1$ and $p_2$ and assuming that $p_1$ differs from $p_2$ on a dense subset of $I$, a necessary and sufficient condition for the equality above is that the unknown functions are of the form
$$
	f=\frac{u}{v},\qquad g=\frac{w}{z},\qquad \mbox{and}\qquad
	p_1q_1=p_2q_2=vz,
$$
where $u,v,w,z:I\to\mathbb{R}$ are arbitrary solutions of the second-order linear differential equation $F''=\gamma F$ ($\gamma\in\R$ is arbitrarily fixed) such that $v>0$ and $z>0$ holds on $I$ and $\{u,v\}$ and $\{w,z\}$ are linearly independent.
\end{abstract}

\maketitle

\section{Introduction} 

Throughout this paper, the symbols $\R$ and $\R_+$ will stand for the sets of real and positive real numbers, respectively, and $I$ will always denote a nonempty open real interval.

Given a strictly monotone continuous function $f:I\to\R$ and a pair of positive valued functions $p=(p_1,p_2):I\to\R_{+}^2$, the \emph{$2$-variable generalized Bajraktarevi\'c mean} $A_{f,p}:I^2\to I$ is given by the following formula:
\Eq{BM}{
	A_{f,p}(x,y):=f^{-1}\bigg(\frac{p_1(x)f(x)+p_2(y)f(y)}{p_1(x)+p_2(y)}\bigg) \qquad (x,y\in I).
}
This is an extension of the notion introduced by Bajraktarevi\'c in \cite{Baj58} and \cite{Baj63}.
It is easy to see that $A_{f,p}$ is a strict mean, i.e., 
\Eq{*}{
	\min(x,y)\leq A_{f,p}(x,y)\leq \max(x,y)\qquad (x,y\in I)
}
holds and the inequalities are strict if $x\neq y$.

The equality problem of generalized Bajraktarevi\'c means has been solved by the authors in \cite{GruPal20}. The main goal of this paper is to investigate the invariance equation of two $2$-variable generalized Bajraktarevi\'c means with respect to the arithmetic mean, i.e., to solve the functional equation
\Eq{inv}{
	A_{f,p}(x,y)+A_{g,q}(x,y)=x+y\qquad (x,y\in I),
}
where $f,g:I\to\R$ are strictly monotone and continuous functions and $p,q:I\to\R_+^2$.

In general classes of means the invariance equation was studied by numerous authors in several papers. Dar\'oczy and P\'ales completely solved this equation for the case of Hölder means in \cite{DarPal02c}. The more general invariance equation for quasi-arithmetic means was first solved by Sut\^o in \cite{Sut14a}, \cite{Sut14b} assuming infinitely many times differentiability and later by Matkowski in \cite{Mat99a} supposing twice continuous differentiability. Finally Dar\'oczy and P\'ales solved this problem without any unnecessary differentiability assumptions in \cite{DarPal02c}. Jarczyk and Matkowski in \cite{JarMat06} and Burai in \cite{Bur07a} has studied the invariance equation in case of three weighted arithmetic means. Jarczyk in \cite{Jar07} gave a solution to this problem with no additional regularity assumptions. Matkowski in \cite{Mat04c} studied the invariance of Cauchy means with respect to the arithmetic mean. The invariance of the arithmetic mean with respect to generalized quasi-arithmetic means (which were introduced by Matkowski in \cite{Mat10b}) was described in \cite{BajPal09a} by Baj\'ak and P\'ales. They also solved the invariance equations of 2-variable Gini and Stolarsky means in \cite{BajPal09b} and \cite{BajPal10}, respectively.

The main difficulty in solving \eq{inv} is that we had to compute the partial derivatives of the two means at the diagonal points of $I^2$ up to the order four and then to solve the system of differential equations so obtained.

In the subsequent section we prove a lemma which allows us to weaken the regularity assumptions in our main theorem. After that the definition and fundamental properties of the Schwarzian derivative are recalled. Then we give a sufficient condition for the invariance equation to be valid. This statement is proved without any regularity assumption on the unknown functions. Finally, assuming four times continuous differentiability, we prove that the same condition is also necessary if $p_1=p_2$ is valid only on a nowhere dense subset of $I$. The case when $p_1=p_2$ holds in a subinterval of $I$ remains open.

\section{Auxiliary results and notation}

Define the \emph{diagonal} $\diag(I^2)$ of $I^2$ and the map $\Delta_2:I\to\diag(I^2)$ by 
\Eq{*}{
\diag(I^2):=\{(x,x)\in\R  ^2\mid x\in I\}\qquad\mbox{and}\qquad
  \Delta_2(x):=(x,x) \qquad(x\in I),
}
and denote $e_1:=(1,0)$ and $e_2:=(0,1)$. Given $p=(p_1,p_2):I\to\R_{+}^2$ and $q=(q_1,q_2):I\to\R_{+}^2$, we will also use the following notations:
\Eq{*}{
	p_0:=p_1+p_2,\qquad q_0:=q_1+q_2,\qquad\mbox{and}\qquad r_0:=\frac{q_0}{p_0}.
} 
We say that $M:I^2\to I$ is a \emph{strict mean} on $I$ if $\min(x,y)<M(x,y)<\max(x,y)$ hold for $x,y\in I$ with $x\neq y$. Furthermore, we call $M$ separately continuous, if the function $x\mapsto M(x,y)$ is continuous for every $y\in I$ and the function $y\mapsto M(x,y)$ is continuous for every $x\in I$. 

With the aid of the following lemma, we can reduce the regularity assumptions in our statements. 
\Lem{reg}{
        Let $M:I^2\to I$ be a strict mean on $I$, let $f:I\to\R$ be a continuous strictly monotone function and $p=(p_1,p_2):I\to\R_{+}^2$. Assume that there exists an open set $U\subseteq I^2$ containing $\diag(I^2)$ such that $A_{f,p}=M$ holds on $U$.  Then the following two assertions hold.
        \begin{enumerate}[(i)]
         \item If $M$ is separately continuous on $U$, then $p_1$ and $p_2$ are continuous on $I$.
         \item Let $k\in\N$. Assume that $M$ is $k$ times partially differentiable (resp.\ $k$ times continuously partially differentiable) with respect to its first and second variables on $U$ and $f:I\to\R$ is a $k$ times differentiable (resp.\ $k$ times continuously differentiable) function on $I$. Then $p_1$ and $p_2$ are $k$ times differentiable (resp.\ $k$ times  continuously differentiable) on $I$.
        \end{enumerate}
}

\begin{proof} By our assumption, for all $(x,y)\in U$, we have that
\Eq{*}{
  A_{f,p}(x,y)=M(x,y).
}
This is equivalent to the following equality
\Eq{ee}{
   \frac{p_1(x)f(x)+p_2(y)f(y)}{p_1(x)+p_2(y)}
   =f(M(x,y)) \qquad((x,y)\in U).
}
Observe that, for $x,y\in I$ with $x\neq y$, the inequalities $\min(x,y)<M(x,y)<\max(x,y)$ and the strict monotonicity of $f$ imply that $f(y)\neq f(M(x,y))$. Thus, solving equation \eq{ee} with respect to $p_2(y)$, we get
\Eq{qi}{
  p_2(y)=p_1(x)\frac{f(M(x,y))-f(x)}{f(y)-f(M(x,y))}
  \qquad((x,y)\in U,\,x\neq y).
}
Let $x_0\in I$ be an arbitrarily fixed point. The pair $(x_0,x_0)$ is an interior point of $U$, therefore, there exists $x\in I\setminus\{x_0\}$ such that $(x,x_0)\in U$. Then the set 
\Eq{*}{
  V:=\{y\in I\mid (x,y)\in U,\,x\neq y\} 
}
is a neighborhood of $x_0$ on which we have the equality \eq{qi} for $p_2$. On the other hand, the continuity of $M$ in its second variable implies that the right hand side of \eq{qi} is a continuous function of $y$ on $V$. Therefore, $p_2$ is continuous at $x_0$, resulting that $p_2$ is continuous on $I$. A similar argument shows that $p_1$ is continuous due to the continuity of $M$ in its first variable.

By the standard calculus rules, under the $k$-times (continuous) differentiability assumptions, the right hand side of \eq{qi} is also $k$-times (continuously) differentiable on the above constructed set $V$, in particular, at $x_0$, which yields the same property for $p_2$ at $x_0$. The arbitrariness of $x_0$ shows that $p_2$ is $k$-times (continuously) differentiable on $I$.
\end{proof}

\section{The Schwarzian derivative}

For a three times differentiable function $f: I\to\R$ with a nonvanishing first derivative, we recall the notion of the Schwarzian derivative $S(f): I\to\R$ which is defined by the following formula:
\Eq{*}{
	S(f)=\frac{f'''}{f'}-\frac{3}{2}\bigg(\frac{f''}{f'}\bigg)^2.
}

\Lem{uv}{Let $\gamma\in\R$ and let $f:I\to\R$ be a three times differentiable function such that $f'$ does not vanish on $I$. Then the following two assertions are equivalent. 
\begin{enumerate}[(i)]
 \item There exist twice differentiable functions $u,v:I\to\R$ such that
 $v$ does not vanish on $I$,
 \Eq{*}{
   u''=\gamma u,\qquad v''=\gamma v,\qquad\mbox{and}\qquad f=\frac{u}{v}.
 }
 \item $f$ satisfies the third-order differential equation
 \Eq{*}{
   S(f)=-2\gamma.
 }
\end{enumerate}
}

\begin{proof}
Assume that assertion (i) holds. Then $u$ and $v$ are three times differentiable functions and
\Eq{*}{
u&=fv, \\ u'&=f'v+fv',\\ 
\gamma u= u''&=f''v+2f'v'+fv''=(\gamma f+f'')v+2f'v',\\
\gamma u'= u'''&=f'''v+3f''v'+3f'v''+fv'''
=(3\gamma f'+f''')v+(\gamma f+3f'')v'.
}
This is a system of homogeneous linear equations with respect to the unknowns $(u,u',v,v')$. Thus its base determinant has to be zero, that is,
\Eq{*}{
0&=\left|\begin{matrix}
 1 & 0 & f & 0 \\
 0 & 1 & f' & f \\
 \gamma & 0 & \gamma f+f''& 2f' \\
 0 & \gamma & 3\gamma f'+f''' & \gamma f+3f''
\end{matrix}\right|
=\left|\begin{matrix}
 1 & 0 & 0 & 0 \\
 0 & 1 & f' & 0 \\
 \gamma & 0 & f''& 2f' \\
 0 & \gamma & 3\gamma f'+f''' & 3f''
\end{matrix}\right|
=\left|\begin{matrix}
  1 & f' & 0 \\
  0 & f''& 2f' \\
  \gamma & 3\gamma f'+f''' & 3f''
\end{matrix}\right|\\
&=\left|\begin{matrix}
  1 & 0 & 0 \\
  0 & f''& 2f' \\
  \gamma & 2\gamma f'+f''' & 3f''
\end{matrix}\right|
=\left|\begin{matrix}
f''& 2f' \\
2\gamma f'+f''' & 3f''
\end{matrix}\right|
=3f''^2-4\gamma f'^2-2f'f'''
=-2f'^2(S(f)+2\gamma).
}
Therefore, assertion (ii) is valid.

To prove the converse, assume that assertion (ii) holds. Let $x_0\in I$ be fixed and let $u$ and $v$ be linearly independent functions that satisfy $u''=\gamma u$, $v''=\gamma v$ and
\Eq{*}{
   u(x_0)&=f(x_0),&\qquad   
   v(x_0)&=1,\\
   u'(x_0)&=f'(x_0)-\frac12f'(x_0)^{-1}f''(x_0)f(x_0),&\qquad   
   v'(x_0)&=-\frac12f'(x_0)^{-1}f''(x_0).
}
By the Liouville Theorem, the Wronskian $W(u,v):=u'v-v'u$ of the linear differential equation $F''=\gamma F$ is constant. Therefore, $W(u,v)\equiv W(u,v)(x_0)=f'(x_0)$. Let $V$ be the largest open subinterval of $I$ containing $x_0$ such that $v$ is positive on $V$. Define $g:V\to\R$ by $g:=u/v$. We can see that
\Eq{*}{
  g(x_0)&=\frac{u(x_0)}{v(x_0)}=f(x_0),\\
  g'(x_0)&=\Big(\frac{u'v-v'u}{v^2}\Big)(x_0)=f'(x_0),\\
  g''(x_0)&=\Big(\frac{u'v-v'u}{v^2}\Big)'(x_0)
  =\bigg(\frac{f'(x_0)}{v^2}\bigg)'(x_0)
  =-2f'(x_0)\Big(\frac{v'}{v^2}\Big)(x_0)=f''(x_0).
}
On the other hand, by the first part of the proof, $g$ also satisfies the differential equation $S(g)=-2\gamma$ on $V$. Thus $f$ and $g$ are solutions of the same ordinary differential equation and they satisfy the same initial value condition at $x_0$. By the existence and uniqueness theorem for ordinary differential equations it follows that $f=g$ on $V$. If $V$ were a proper subinterval of $I$, then one of its endpoints, say the lower one, would belong to $I$. At this endpoint, the function $v$ vanishes, hence the right limit of $g$ does not exist contradicting that the right limit of $f$ exists at this point. This contradiction shows that $V=I$ and hence, $f$ is of the form stated in assertion (i).
\end{proof}

The next result is the particular case of \lem{uv} when $\gamma=0$.

\Cor{S(f)=0}{Let $f:I\to\R$ be a three times differentiable function such that $f'$ does not vanish on $I$. Then
\Eq{*}{
 S(f)=0
}
holds on $I$ if and only if there exist four constants $a,b,c,d\in\R$ with $ad\neq bc$ and $0\not\in cI+d$ such that 
\Eq{*}{
  f(x)=\frac{ax+b}{cx+d} \qquad(x\in I).
}}

\section{A sufficient condition for the invariance equation}

In what follows we prove a sufficient condition for $(f,p),(g,q)$ to be a solution of \eq{inv}. 

For a real number $\gamma\in\R$, we introduce the sine and cosine type functions $S_\gamma,C_\gamma:\R\to\R$ by
\Eq{*}{
  S_\gamma(x)&:=\sum_{k=0}^\infty\frac{\gamma^kx^{2k+1}}{(2k+1)!}=\begin{cases}
           \dfrac{\sin(\sqrt{-\gamma}x)}{\sqrt{-\gamma}} 
           & \mbox{ if } \gamma<0, \\
           x & \mbox{ if } \gamma=0, \\
           \dfrac{\sinh(\sqrt{\gamma}x)}{\sqrt{\gamma}} 
           & \mbox{ if } \gamma>0, \\
         \end{cases}\\
  C_\gamma(x)&:=\sum_{k=0}^\infty\frac{\gamma^kx^{2k}}{(2k)!}
   =\begin{cases}
           \cos(\sqrt{-\gamma}x) & \mbox{ if } \gamma<0, \\
           1 & \mbox{ if } \gamma=0, \\
           \cosh(\sqrt{\gamma}x) & \mbox{ if } \gamma>0. \\
        \end{cases}
}

\Thm{Suff}{Let $\gamma\in\R$ be a real constant, let $u,v,w,z:I\to\R$ be arbitrary solutions of the second-order linear differential equation $F''=\gamma F$ such that $v>0$ and $z>0$ holds on $I$ and $\{u,v\}$ and
$\{w,z\}$ are linearly independent. Assume that the functions $f,g:I\to\R$, $p=(p_1,p_2):I\to\R_{+}^2$, and $q=(q_1,q_2):I\to\R_{+}^2$ satisfy
\Eq{fgpq}{
  f=\frac{u}{v},\qquad g=\frac{w}{z},\qquad \mbox{and}\qquad
  p_1q_1=p_2q_2=vz. 
}
Then $f$ and $g$ are strictly monotone and continuous and the invariance equation \eq{inv} holds for all $x,y\in I$.}

\begin{proof} 
By basic results on linear homogeneous differential equations, the functions $S_\gamma$ and $C_\gamma$ form a fundamental system of solutions for the differential equation $F''=\gamma F$. Therefore, the pairs $(u,v)$ and $(w,z)$ are equivalent to $(S_\gamma,C_\gamma)$, that is, there exist $a_1,b_1,c_1,d_1,a_2,b_2,c_2,d_2\in\R$ real constants such that $a_1d_1\neq b_1c_1$, $a_2d_2\neq b_2c_2$ and
\Eq{*}{
    u=a_1 S_\gamma+b_1 C_\gamma,\qquad 
    v=c_1 S_\gamma+d_1 C_\gamma, \qquad\mbox{and}\qquad 
    w=a_2 S_\gamma+b_2 C_\gamma,\qquad 
    z=c_2 S_\gamma+d_2 C_\gamma.
}
In view of the sufficiency part of \cite[Theorem 6]{PalZak19}, this implies the identity
\Eq{uvwz}{
  \Big(\frac{u}{v}\Big)^{-1}\bigg(\frac{tu(x)+su(y)}{tv(x)+sv(y)}\bigg)
  +\Big(\frac{w}{z}\Big)^{-1}\bigg(\frac{sw(x)+tw(y)}{sz(x)+tz(y)}\bigg)=x+y
  \qquad(x,y\in I,\,t,s>0).
}
On the other hand, with $t:=\frac{p_1(x)}{v(x)}=\frac{z(x)}{q_1(x)}$ and $s:=\frac{p_2(y)}{v(y)}=\frac{z(y)}{q_2(y)}$, using \eq{fgpq}, we have
\Eq{*}{
  \Big(\frac{u}{v}\Big)^{-1}\bigg(\frac{tu(x)+su(y)}{tv(x)+sv(y)}\bigg)
  =f^{-1}\bigg(\frac{p_1(x)f(x)+p_2(y)f(y)}{p_1(x)+p_2(y)}\bigg)
  =A_{f,p}(x,y)
}
and
\Eq{*}{
  \Big(\frac{w}{z}\Big)^{-1}\bigg(\frac{sw(x)+tw(y)}{sz(x)+tz(y)}\bigg)
  =g^{-1}\bigg(\frac{q_1(x)g(x)+q_2(y)g(y)}{q_1(x)+q_2(y)}\bigg)
  =A_{g,q}(x,y).
}
Therefore, \eq{uvwz} yields that \eq{inv} is satisfied.
\end{proof}

\section{Partial derivatives of generalized Bajraktarevi\'c means}

In the next result we recall the formulas for the partial derivatives of generalized Bajraktarevi\'c means up to third-order at diagonal points of $I^2$. These assertions were proved in \cite{GruPal20} under tight regularity assumptions. We also calculate the fourth-order partial derivative $\partial_1^2 \partial_2^2 A_{f,p}$.

\Thm{DB}{
	Let $\ell\in\{1,2,3,4\}$, let $f: I\to\R$ be an $\ell$ times differentiable function on $I$ with a nonvanishing first derivative, and let $p=(p_1,p_2):I\to\R_{+}^2$. Then we have the following assertions.
	\begin{enumerate}[(i)]
		\item If $\ell=1$, $i\in \{1,2\}$, and $p_i$ is continuous on $I$, then the first-order partial derivative $\partial_i A_{f,p}$ exists on $\diag(I^2)$ and
		\Eq{*}{
		\partial_i A_{f,p}\circ\Delta_2=\frac{p_i}{p_0}.
		}
		\item If $\ell=2$, $p_1$ and $p_2$ are differentiable on $I$, then the second-order partial derivative $\partial_1\partial_2 A_{f,p}$ exists on $\diag(I^2)$ and
		\Eq{*}{
		\partial_1\partial_2 A_{f,p}\circ\Delta_2 
		=-\frac{p_1p_2}{p_0^2}\bigg(\frac{(p_1p_2)'}{p_1p_2}+\frac{f''}{f'}\bigg).
		}
		\item If $\ell=3$ and $p_1$ and $p_2$ are twice differentiable on $I$, then the third-order partial derivative $\partial_1^2\partial_2 A_{f,p}$ exists on $I^2$ and
		\Eq{*}{
		\partial_1^2 \partial_2 A_{f,p}\circ\Delta_2
		&=-\frac14\Big(\frac{p_1-p_2}{p_0}\Big)''-\frac{3p_0(p_1-p_2)}{16p_1p_2}\bigg(\Big(\frac{p_1-p_2}{p_0}\Big)'\bigg)^2-\frac12\bigg(\frac{p_1p_2}{p_0^2}\bigg(\frac{(p_1p_2)'}{p_1p_2}+\frac{f''}{f'}\bigg)\bigg)'\\
        &\qquad +\frac{3p_1p_2}{4p_0^2}\cdot\frac{p_1-p_2}{p_0}\bigg(\frac{(p_1p_2)'}{p_1p_2}+\frac{f''}{f'}\bigg)^2-\frac{p_1p_2}{2p_0^2}\cdot\frac{p_1-p_2}{p_0}S(f).
        }
        \item If $\ell=4$, $p_1$ and $p_2$ are twice continuously differentiable on $I$, then the fourth-order partial derivative $\partial_1^2\partial_2^2 A_{f,p}$ exists on $\diag(I^2)$ and
		\Eq{*}{
		\partial_1^2\partial_2^2 A_{f,p}\circ\Delta_2
		&=\bigg(\Big(\frac{p_1p_2}{p_0^2}\Big)''
        +\frac38\Big(6-\frac{p_0^2}{p_1p_2}\Big)\bigg(\Big(\frac{p_1-p_2}{p_0}\Big)'\bigg)^2\bigg)\bigg(\frac{(p_1p_2)'}{p_1p_2}+\frac{f''}{f'}\bigg)\\
        &\qquad-\Big(\frac{p_1p_2}{p_0^2}\Big)'\bigg(\frac{(p_1p_2)'}{p_1p_2}+\frac{f''}{f'}\bigg)'
        -\frac12\Big(\frac{p_1p_2}{p_0^2}\Big)^2\bigg(6-\frac{p_0^2}{p_1p_2}\bigg)\bigg(\frac{(p_1p_2)'}{p_1p_2}+\frac{f''}{f'}\bigg)^3\\
        &\qquad+\Big(\frac{p_1p_2}{p_0^2}\Big)^2\bigg(6-\frac{p_0^2}{p_1p_2}\bigg)\bigg(\frac{(p_1p_2)'}{p_1p_2}+\frac{f''}{f'}\bigg)S(f)-\Big(\frac{p_1p_2}{p_0^2}\Big)^2 S(f)'.
        }
	\end{enumerate}
}

\begin{proof} It follows from elementary calculus rules that partial derivative $\partial_1^\alpha\partial_2^\beta A_{f,p}$ exists on $I^2$ if $f$ is $(\alpha+\beta)$-times differentiable with a nonvanishing first derivative, furthermore, $p_1$ is $\alpha$-times and $p_2$ is $\beta$-times differentiable on $I$. 

By equality \eq{ee}, with the notation $M:=A_{f,p}$, for all $x,y\in I$, we have
\Eq{*}{
  p_1(x)f(x)+p_2(y)f(y)
  =p_1(x)f\big(M(x,y)\big)+p_2(y)f\big(M(x,y)\big).
}
Let $\ell\in\{1,2,3,4\}$, $\alpha,\beta\in\{0,1,2\}$ and let $\delta_{\cdot,\cdot}$ denote the standard Kronecker symbol. Differentiating this equality with respect to the first variable $\alpha$-times and the second variable $\beta$-times by applying the generalized Leibniz Product Rule, we get
\Eq{*}{
   \delta_{0,\beta}\sum_{i=0}^\alpha &\binom{\alpha}{i}
   p_1^{(\alpha-i)}(x)f^{(i)}(x)
   +\delta_{0,\alpha}\sum_{j=0}^\beta \binom{\beta}{j}
   p_2^{(\beta-j)}(y)f^{(j)}(y)\\
   &=\sum_{i=0}^{\alpha}\binom{\alpha}{i}
   p_1^{(\alpha-i)}(x)
   \cdot\partial_1^i\partial_2^\beta (f\circ M)(x,y)
   +\sum_{j=0}^{\beta}\binom{\beta}{j}p_2^{(\beta-j)}(y)
   \cdot\partial_1^\alpha\partial_2^j (f\circ M)(x,y).
}
Restricting this equality to the diagonal of $I^2$, we obtain
\Eq{LPR}{
   \delta_{0,\beta}\sum_{i=0}^\alpha &\binom{\alpha}{i}
   p_1^{(\alpha-i)}f^{(i)}
   +\delta_{0,\alpha}\sum_{j=0}^\beta \binom{\beta}{j}
   p_2^{(\beta-j)}f^{(j)}\\
   &=\sum_{i=0}^{\alpha}\binom{\alpha}{i}
   p_1^{(\alpha-i)}
   \cdot\partial_1^i\partial_2^\beta (f\circ M)\circ\Delta_2
   +\sum_{j=0}^{\beta}\binom{\beta}{j}p_2^{(\beta-j)}
   \cdot\partial_1^\alpha\partial_2^j (f\circ M)\circ\Delta_2.
}
For the computation of the partial derivatives $\partial_1^i\partial_2^j (f\circ M)$, the following easy-to-see formulas apply:
\Eq{PDS}{
  \partial_\mu (f\circ M)&=(f'\circ M)\cdot\partial_\mu M, \\
  \partial_\mu\partial_\nu (f\circ M)
  &=(f''\circ M)\cdot\partial_\mu M\cdot\partial_\nu M+(f'\circ M)\cdot\partial_\mu\partial_\nu M,\\
  \partial_\mu^2\partial_\nu(f\circ M)
  &=(f'''\circ M)\cdot\partial_\mu M^2 \cdot\partial_\nu M+(f''\circ M)\cdot\big(\partial_\mu^2 M\cdot\partial_\nu M+2\partial_\mu\partial_\nu M\cdot\partial_\mu M\big) \\
  &\quad+(f'\circ M)\cdot\partial_\mu^2\partial_\nu M,\\
  \partial_1^2\partial_2^2 (f\circ M)
  &=(f''''\circ M)\cdot\partial_1 M^2\cdot\partial_2 M^2+(f'\circ M)\cdot\partial_1^2 \partial_2^2 M \\
  &\quad+(f'''\circ M)\cdot\big(\partial_1^2 M\cdot\partial_2 M^2+\partial_2^2 M\cdot\partial_1 M^2+4\partial_1\partial_2 M\cdot\partial_1 M\cdot\partial_2 M \big) \\
  &\quad+(f''\circ M)\cdot\big(\partial_1^2 M\cdot\partial_2^2 M+2\partial_1\partial_2 M^2+2\partial_1^2\partial_2 M\cdot\partial_2 M+2\partial_1\partial_2^2 M\cdot\partial_1 M\big),
}
where $\mu,\nu\in\{1,2\}$.

In the particular case when $\alpha=1$ and $\beta=0$, \eq{LPR} and also the first formula from \eq{PDS} yield
\Eq{*}{
  p'_1f+p_1f'=p'_1f+p_0f'\cdot\partial_1M\circ\Delta_2.
}
Expressing $\partial_1M\circ\Delta_2$, we arrive at assertion (i) for $i=1$. The case $i=2$ can be seen similarly.

In the case $\alpha=\beta=1$, using the first two formulas from \eq{PDS}, the equality \eq{LPR} simplifies to
\Eq{*}{
  0=p'_1f'\cdot\partial_2M\circ\Delta_2+p'_2f'\cdot\partial_1M\circ\Delta_2
  +p_0\big(f''\cdot(\partial_1M\cdot\partial_2M)\circ\Delta_2+f'\cdot\partial_1\partial_2M\circ\Delta_2\big),
}
which, using (i), implies that
\Eq{d12}{
  f''\cdot(\partial_1 M\cdot\partial_2 M)\circ\Delta_2+f'\cdot\partial_1\partial_2 M\circ\Delta_2
  =-\frac{(p_1p_2)'}{p_0^2}f'.
}
Now, applying the formulas from (i) again, we can conclude that (ii) is valid. Using (i) and (ii), we can also get the formula
\Eq{i2}{
\partial_i^2M\circ\Delta_2
&=\big(\partial_iM(x,x)\circ\Delta_2\big)'-\partial_1\partial_2M\circ\Delta_2\\
&=\Big(\frac{p_i}{p_0}\Big)'+\frac{p_1p_2}{p_0^2}\bigg(\frac{(p_1p_2)'}{p_1p_2}+\frac{f''}{f'}\bigg)
=2\frac{p'_ip^{\phantom1}_{3-i}}{p_0^2}+\frac{p_1p_2}{p_0^2}\cdot\frac{f''}{f'}.
}
Using assertion (i) and \eq{i2}, it follows that
\Eq{dii}{
  f''\cdot(\partial_i M^2)\circ\Delta_2+f'\cdot\partial_i^2 M\circ\Delta_2 =\frac{p_i}{p_0}f''+2\dfrac{p_i' p^{\phantom 1}_{3-i}}{p_0^2}f'.
}

In the case when $\alpha=2$, $\beta=1$, the equality \eq{LPR} yields
\Eq{*}{
   0=p_1''\cdot\partial_2 (f\circ M)\circ\Delta_2
   +2p_1'\cdot\partial_1\partial_2 (f\circ M)\circ\Delta_2
   +p_2'\cdot\partial_1^2(f\circ M)\circ\Delta_2
   +p_0\cdot\partial_1^2\partial_2 (f\circ M)\circ\Delta_2.
}
In view of \eq{PDS}, we can rewrite this equality as
\Eq{*}{
   0&=p_1''f'\cdot\partial_2 M\circ\Delta_2+2p_1'\Big(f''\cdot(\partial_1 M\cdot\partial_2 M)\circ\Delta_2+f'\cdot\partial_1\partial_2 M\circ\Delta_2\Big)
   \\
   &\qquad+p_2'\Big(f''\cdot(\partial_1 M^2)\circ\Delta_2+f'\cdot\partial_1^2 M\circ\Delta_2\Big)+p_0\Big(f'''\cdot(\partial_1 M^2 \cdot\partial_2 M)\circ\Delta_2\\
   &\qquad+f''\cdot\big(\partial_1^2 M\cdot\partial_2 M+2\partial_1\partial_2 M\cdot\partial_1 M\big)\circ\Delta_2
   +f'\cdot\partial_1^2\partial_2 M\circ\Delta_2\Big).
}
This equality, using assertion (i), \eq{d12} and \eq{dii}, implies that
\Eq{3rdO}{
	f'''\cdot&(\partial_1M^2\cdot\partial_2M)\circ\Delta_2+f''\cdot(\partial_1^2M\cdot\partial_2M
	+2\partial_1\partial_2M\cdot\partial_1M)\circ\Delta_2+f'\cdot\partial_1^2\partial_2M\circ\Delta_2\\
	&=-\frac{p_1''p_2}{p_0^2}f'+2\frac{p_1'(p_1p_2)'}{p_0^3}f'-\frac{p_2'}{p_0}\bigg(\frac{p_1}{p_0}f''+2\dfrac{p_1' p^{\phantom 1}_2}{p_0^2}f'\bigg)\\
   &=\frac{2p_1'((p_1p_2)'-p_2' p^{\phantom 1}_2)-p_1''p_2p_0}{p_0^3}f'-\frac{p_2'p_1}{p_0^2}f''.
}
Using assertions (i), (ii), we can also get
\Eq{112}{
(\partial_1^2M\cdot\partial_2M+2\partial_1\partial_2M\cdot\partial_1M)\circ\Delta_2  =2\frac{p'_1p^2_{2}-p_1(p_1p_2)'}{p_0^3}+\frac{p_1p_2(p_2-2p_1)}{p_0^3}\cdot\frac{f''}{f'}.
}
Thus, dividing equation \eq{3rdO} by $f'$ side by side, it reduces to
\Eq{*}{
	\frac{p_1^2p_2}{p_0^3}\cdot\frac{f'''}{f'}+\bigg(2\frac{p'_1p^2_{2}-p_1(p_1p_2)'}{p_0^3}&+\frac{p_1p_2(p_2-2p_1)}{p_0^3}\cdot\frac{f''}{f'}\bigg)\frac{f''}{f'}+\partial_1^2\partial_2M\circ\Delta_2\\
   &=\frac{2p_1'((p_1p_2)'-p_2' p^{\phantom 1}_2)-p_1''p_2p_0}{p_0^3}-\frac{p_2'p_1}{p_0^2}\cdot\frac{f''}{f'}.
}
Therefore,
\Eq{*}{
  \partial_1^2\partial_2M\circ\Delta_2
  &=\frac{2p_1'((p_1p_2)'-p_2' p^{\phantom 1}_2)-p_1''p_2p_0}{p_0^3}-\frac{(2p'_1p_2+p_2'p_1)(p_2-p_1)}{p_0^3}\cdot\frac{f''}{f'}\\
  &\qquad -\frac{p_1p_2(p_2-2p_1)}{p_0^3}\bigg(\frac{f''}{f'}\bigg)^2
  -\frac{p_1^2p_2}{p_0^3}\cdot\frac{f'''}{f'}.
}
From here, assertion (iii) follows.

In the case when $\alpha=\beta=2$, the equality \eq{LPR} yields
\Eq{*}{
  0=&p_2''\cdot\partial_1^2 (f\circ M)\circ\Delta_2
   +p_1''\cdot\partial_2^2 (f\circ M)\circ\Delta_2\\
   &+2p_2'\cdot\partial_1^2\partial_2(f\circ M)\circ\Delta_2
   +2p_1'\cdot\partial_1\partial_2^2 (f\circ M)\circ\Delta_2
   +p_0\cdot\partial_1^2\partial_2^2 (f\circ M)\circ\Delta_2.
}
Hence 
\Eq{*}{
	0&=p_2''\Big(f''\cdot(\partial_1M^2)\circ\Delta_2+f'\cdot\partial_1^2M\circ\Delta_2\Big)+p_1''\Big(f''\cdot(\partial_2M^2)\circ\Delta_2+f'\cdot\partial_2^2M\circ\Delta_2\Big)\\ &\quad+2p_2'\Big(f'''\cdot(\partial_1M^2\cdot\partial_2M)\circ\Delta_2+f''\cdot(\partial_1^2M\cdot\partial_2M
	+2\partial_1\partial_2M\cdot\partial_1M)\circ\Delta_2+f'\cdot\partial_1^2\partial_2M\circ\Delta_2\Big)\\
	&\quad+2p_1'\Big(f'''\cdot(\partial_2M^2\cdot\partial_1M)\circ\Delta_2+f''\cdot(\partial_2^2M\cdot\partial_1M+2\partial_1\partial_2M\cdot\partial_2M)\circ\Delta_2+f'\cdot\partial_2^2\partial_1M\circ\Delta_2\Big)\\
	&\quad+p_0\Big(f''''\cdot(\partial_1 M^2\cdot\partial_2 M^2)\circ\Delta_2+f'\cdot\partial_1^2 \partial_2^2 M\circ\Delta_2 \\
	&\qquad\qquad+f'''\cdot\big(\partial_1^2 M\cdot\partial_2 M^2+\partial_2^2 M\cdot\partial_1 M^2+4\partial_1\partial_2 M\cdot\partial_1 M\cdot\partial_2 M \big)\circ\Delta_2 \\
	&\qquad\qquad+f''\cdot\big(\partial_1^2 M\cdot\partial_2^2 M+2\partial_1\partial_2 M^2+2\partial_1^2\partial_2 M\cdot\partial_2 M+2\partial_1\partial_2^2 M\cdot\partial_1 M\big)\circ\Delta_2\Big).
}
Using \eq{112} and its symmetric counterpart, we can obtain that
\Eq{*}{
  &\big(\partial_1^2 M\cdot\partial_2 M^2+\partial_2^2 M\cdot\partial_1 M^2+4\partial_1\partial_2 M\cdot\partial_1 M\cdot\partial_2 M \big)\circ\Delta_2\\
  &=\big((\partial_1^2 M\cdot\partial_2 M+2\partial_1\partial_2 M\cdot\partial_1 M)\cdot\partial_2 M\big)\circ\Delta_2+\big((\partial_2^2 M\cdot\partial_1 M+2\partial_1\partial_2 M\cdot\partial_2 M)\cdot\partial_1 M\big)\circ\Delta_2\\
  &=\bigg(2\frac{p'_1p^2_{2}-p_1(p_1p_2)'}{p_0^3}+\frac{p_1p_2(p_2-2p_1)}{p_0^3}\Cdot\frac{f''}{f'}\bigg)\frac{p_2}{p_0}+\bigg(2\frac{p'_2p^2_1-p_2(p_1p_2)'}{p_0^3}+\frac{p_1p_2(p_1-2p_2)}{p_0^3}\Cdot\frac{f''}{f'}\bigg)\frac{p_1}{p_0}\\
  &=2\frac{p_1'p_2^3+p_2'p_1^3-2p_1p_2(p_1p_2)'}{p_0^4}+\frac{p_1p_2(p_0^2-6p_1p_2)}{p_0^4}\cdot\frac{f''}{f'}
}
and
\Eq{*}{
  &\big(\partial_1^2 M\cdot\partial_2^2 M+2\partial_1\partial_2 M^2+2\partial_1^2\partial_2 M\cdot\partial_2 M+2\partial_1\partial_2^2 M\cdot\partial_1 M\big)\circ\Delta_2\\
  &=\bigg(2\frac{p'_1p^{\phantom1}_{2}}{p_0^2}+\frac{p_1p_2}{p_0^2}\cdot\frac{f''}{f'}\bigg)\bigg(2\frac{p'_2p^{\phantom1}_{1}}{p_0^2}+\frac{p_1p_2}{p_0^2}\cdot\frac{f''}{f'}\bigg)+2\bigg(\frac{(p_1p_2)'}{p_0^2}+\frac{p_1p_2}{p_0^2}\cdot\frac{f''}{f'}\bigg)^2 \\
  &\qquad+2\bigg(\frac{2p_1'((p_1p_2)'-p_2' p^{\phantom 1}_2)-p_1''p_2p_0}{p_0^3}-\frac{(2p'_1p_2+p_2'p_1)(p_2-p_1)}{p_0^3}\cdot\frac{f''}{f'} -\frac{p_1p_2(p_2-2p_1)}{p_0^3}\bigg(\frac{f''}{f'}\bigg)^2  \\
  &\qquad\qquad\quad-\frac{p_1^2p_2}{p_0^3}\cdot\frac{f'''}{f'}\bigg)\frac{p_2}{p_0} +2\bigg(\frac{2p_2'((p_1p_2)'-p_1' p^{\phantom 1}_1)-p_2''p_1p_0}{p_0^3}-\frac{(2p'_2p_1+p_1'p_2)(p_1-p_2)}{p_0^3}\cdot\frac{f''}{f'}\\
  &\qquad\qquad\qquad\qquad\qquad\qquad\qquad\quad -\frac{p_1p_2(p_1-2p_2)}{p_0^3}\bigg(\frac{f''}{f'}\bigg)^2 -\frac{p_1p_2^2}{p_0^3}\cdot\frac{f'''}{f'}\bigg)\frac{p_1}{p_0}\\
  &=\frac{4p'_1p_2'(3p_1p_2-p_0^2)+6(p_1p_2)'^2}{p_0^4}
  -2\frac{p_1''p_2^2+p_2''p_1^2}{p_0^3}\\
  &\qquad+\frac{4(p_1p_2)'(5p_1p_2-p_0^2)+2p_1p_2(p_1'p_1+p_2'p_2)}{p_0^4}\Cdot\frac{f''}{f'}
  +\frac{p_1p_2(15p_1p_2-2p_0^2)}{p_0^4}\bigg(\frac{f''}{f'}\bigg)^2-4\frac{p_1^2p_2^2}{p_0^4}\Cdot\frac{f'''}{f'}.
}

Dividing this equality by $p_0f'$ side by side, and then using assertions (i), (ii) and \eq{i2} for the computation of the at most second-order partial derivatives, a simple computation yields that
\Eq{*}{
	0&=\frac{p_2''}{p_0}\bigg(2\dfrac{p_1' p^{\phantom 1}_{2}}{p_0^2}+\frac{p_1}{p_0}\cdot\frac{f''}{f'}\bigg)+2\frac{p_2'}{p_0}\bigg(\frac{2p_1'((p_1p_2)'-p_2' p^{\phantom 1}_2)-p_1''p_2p_0}{p_0^3}-\frac{p_2'p_1}{p_0^2}\cdot\frac{f''}{f'}\bigg)
	\\ 
	&\quad+\frac{p_1''}{p_0}\bigg(2\dfrac{p_2' p^{\phantom 1}_{1}}{p_0^2}+\frac{p_2}{p_0}\cdot\frac{f''}{f'}\bigg)
	+2\frac{p_1'}{p_0}\bigg(\frac{2p_2'((p_1p_2)'-p_1' p^{\phantom 1}_1)-p_2''p_1p_0}{p_0^3}-\frac{p_1'p_2}{p_0^2}\cdot\frac{f''}{f'}\bigg)\\
	&\quad+\frac{p_1^2p_2^2}{p_0^4}\cdot\frac{f''''}{f'}+\bigg(2\frac{p_1'p_2^3+p_2'p_1^3-2p_1p_2(p_1p_2)'}{p_0^4}+\frac{p_1p_2(p_0^2-6p_1p_2)}{p_0^4}\cdot\frac{f''}{f'}\bigg) \frac{f'''}{f'} \\
	&\quad+\bigg(\frac{4p'_1p_2'(3p_1p_2-p_0^2)+6(p_1p_2)'^2}{p_0^4}
  +\frac{4(p_1p_2)'(5p_1p_2-p_0^2)+2p_1p_2(p_1'p_1+p_2'p_2)}{p_0^4}\cdot\frac{f''}{f'}\\
  &\qquad\quad -2\frac{p_1''p_2^2+p_2''p_1^2}{p_0^3}+\frac{p_1p_2(15p_1p_2-2p_0^2)}{p_0^4}\bigg(\frac{f''}{f'}\bigg)^2-4\frac{p_1^2p_2^2}{p_0^4}\cdot\frac{f'''}{f'}\bigg)\frac{f''}{f'}+\partial_1^2 \partial_2^2 M\circ\Delta_2.
}
Therefore,
\Eq{*}{
\partial_1^2 \partial_2^2 & M\circ\Delta_2
=2\frac{(p_1p_2)''p_0'p_0-p_0''(p_1p_2)'p_0-6p_1'p_2'(p_1p_2)'}{p_0^4}\\
	&\quad+\frac{(p_1p_2)''p_0^2-2p_0''p_0p_1p_2+p_1'p_2'(2p_0^2-24p_1p_2)+p_1'^2p_2(2p_0-6p_2)+p_2'^2p_1(2p_0-6p_1)}{p_0^4}\Cdot\frac{f''}{f'}\\
	&\quad+\frac{(p_1p_2)'(4p_0^2-18p_1p_2)-2p_0'p_0p_1p_2}{p_0^4}\bigg(\frac{f''}{f'}\bigg)^2+\frac{p_1p_2(2p_0^2-15p_1p_2)}{p_0^4}\bigg(\frac{f''}{f'}\bigg)^3\\
	&\quad+\frac{(p_1p_2)'(6 p_1p_2-2p_0^2)+2p_0'p_0p_1p_2}{p_0^4} \cdot\frac{f'''}{f'}+\frac{p_1p_2(10p_1p_2-p_0^2)}{p_0^4}\cdot\frac{f''f'''}{(f')^2}-\frac{(p_1p_2)^2}{p_0^4}\cdot\frac{f^{(4)}}{f'}.
}
From here, we can directly conclude that assertion (iv) holds.
\end{proof}

\section{Necessary conditions for the invariance equation}

In the subsequent lemmas we establish the first-, second-, third-, and fourth-order necessary conditions for the validity of the invariance equation \eq{inv}. Finally, we present the main result of our paper in \thm{Main}.

\Lem{1OC}{
Let $f,g: I\to\R$ be differentiable functions on $I$ with nonvanishing first derivatives and $i\in\{1,2\}$. Let $p=(p_1,p_2): I\to\R_{+}^2$ and $q=(q_1,q_2): I\to\R_{+}^2$ be such that $p_i$ and $q_i$ are continuous on $I$. If $\partial_i A_{f,p}+\partial_i A_{g,q}=1$ holds on $\diag(I^2)$, then
\Eq{1OC}{
  \frac{p_1}{p_0}=\frac{q_2}{q_0} \qquad\mbox{and}\qquad
  \frac{p_2}{p_0}=\frac{q_1}{q_0}
}
and hence
\Eq{pq}{
	p_1q_1=p_2q_2
}
holds on $I$.
}

\begin{proof} By formula (i) of \thm{DB}, the equality $(\partial_i A_{f,p}+\partial_i A_{g,q})\circ\Delta_2=1$ can be rewritten as
\Eq{*}{
    \frac{p_i}{p_0}+\frac{q_i}{q_0}=1,
}
which is equivalent to \eq{1OC} and also equivalent to \eq{pq}.
\end{proof}

\Lem{2OC}{
Let $f,g: I\to\R$ be twice differentiable functions on $I$ with nonvanishing first derivatives. Let $p=(p_1,p_2): I\to\R_{+}^2$ and $q=(q_1,q_2): I\to\R_{+}^2$ be differentiable functions on $I$ and assume that \eq{1OC} holds on $I$. If $\partial_1\partial_2 A_{f,p}+\partial_1\partial_2 A_{g,q}=0$ holds on $\diag(I^2)$, then
\Eq{2OC}{
    \frac{(p_1p_2)'}{p_1p_2}+\frac{f''}{f'}
    +\frac{(q_1q_2)'}{q_1q_2}+\frac{g''}{g'}=0.
}
Consequently, there exists a nonzero constant $\delta$ such that
\Eq{ga}{
	p_1q_1=p_2q_2=\sqrt{\frac{\delta}{f'g'}}
}
is valid on $I$.
}

\begin{proof} The equality $(\partial_1\partial_2 A_{f,p}+\partial_1\partial_2 A_{g,q})\circ\Delta_2=0$, in view of formula (ii) of \thm{DB}, is equivalent to
\Eq{*}{
    \frac{p_1p_2}{p_0^2}\bigg(\frac{(p_1p_2)'}{p_1p_2}+\frac{f''}{f'}\bigg)
    +\frac{q_1q_2}{q_0^2}\bigg(\frac{(q_1q_2)'}{q_1q_2}+\frac{g''}{g'}\bigg)=0.
    }
Multiplying this equality by $\frac{p_0^2}{p_1p_2}$, which by \eq{1OC}, is equal to $\frac{q_0^2}{q_1q_2}$, we can easily see that \eq{2OC} holds on $I$. Integrating both sides of the equality \eq{2OC}, we find that there exists a constant $\delta$ such that
\Eq{*}{
  p_1p_2f'q_1q_2g'=\delta.
}
This equality together with \eq{pq} implies that \eq{ga} is also valid.
\end{proof}

\Lem{3OC}{
	Let $f,g: I\to\R$ be three times differentiable functions on $I$ with nonvanishing first derivatives. Let $p=(p_1,p_2): I\to\R_{+}^2$ and $q=(q_1,q_2): I\to\R_{+}^2$ be twice differentiable functions such that \eq{1OC} and \eq{2OC} hold on $I$. If $\partial_1^2\partial_2 A_{f,p}+\partial_1^2\partial_2 A_{g,q}=0$ holds on $\diag(I^2)$, then
	\Eq{3OC}{
	  \frac{p_1-p_2}{p_0}S(f)+\frac{q_1-q_2}{q_0}S(g)=0.
	}
	Consequently,
	\Eq{pS}{
	  (p_1-p_2)(S(f)-S(g))=0
	}
	is valid on $I$.
}

\begin{proof}
In view of \eq{1OC} and \eq{2OC}, we have that
\Eq{ids}{
	\frac{p_1p_2}{p_0^2}=\frac{q_1q_2}{q_0^2},
	\qquad \frac{p_1-p_2}{p_0}=\frac{q_2-q_1}{q_0},\qquad
	\frac{(p_1p_2)'}{p_1p_2}+\frac{f''}{f'}
	=-\bigg(\frac{(q_1q_2)'}{q_1q_2}+\frac{g''}{g'}\bigg).
}
By using formula (iii) of \thm{DB} for $\partial_1^2 \partial_2 A_{f,p}\circ\Delta_2$ and the analogous formula for $\partial_1^2 \partial_2 A_{g,q}\circ\Delta_2$, the equality $\big(\partial_1^2\partial_2 A_{f,p}+\partial_1^2\partial_2 A_{g,q}\big)\circ\Delta_2=0$ can be rewritten as
\Eq{*}{
  0&=-\frac14\Big(\frac{p_1-p_2}{p_0}\Big)''-\frac{3p_0^2}{16p_1p_2}\cdot\frac{p_1-p_2}{p_0}\bigg(\frac{p_1-p_2}{p_0}\bigg)'^2+\frac{3p_1p_2}{4p_0^2}\cdot\frac{p_1-p_2}{p_0}\bigg(\frac{(p_1p_2)'}{p_1p_2}+\frac{f''}{f'}\bigg)^2\\
    &\quad -\frac12\bigg(\frac{p_1p_2}{p_0^2}\bigg(\frac{(p_1p_2)'}{p_1p_2}+\frac{f''}{f'}\bigg)\bigg)'-\frac{p_1p_2}{2p_0^2}\cdot\frac{p_1-p_2}{p_0}S(f)\\
    &\quad-\frac14\Big(\frac{q_1-q_2}{q_0}\Big)''-\frac{3q_0^2}{16q_1q_2}\cdot\frac{q_1-q_2}{q_0}\bigg(\frac{q_1-q_2}{q_0}\bigg)'^2+\frac{3q_1q_2}{4q_0^2}\cdot\frac{q_1-q_2}{q_0}\bigg(\frac{(q_1q_2)'}{q_1q_2}+\frac{g''}{g'}\bigg)^2\\
    &\quad -\frac12\bigg(\frac{q_1q_2}{q_0^2}\bigg(\frac{(q_1q_2)'}{q_1q_2}+\frac{g''}{g'}\bigg)\bigg)'-\frac{q_1q_2}{2q_0^2}\cdot\frac{q_1-q_2}{p_0}S(g).
}
Using the identities in \eq{ids}, this equality is equivalent to 
\Eq{*}{
  0=-\frac{p_1p_2}{2p_0^2}\cdot\frac{p_1-p_2}{p_0}S(f)
   -\frac{q_1q_2}{2q_0^2}\cdot\frac{q_1-q_2}{p_0}S(g).
}
Multiplying the last equation by $\frac{2p_0^2}{p_1p_2}=\frac{2q_0^2}{q_1q_2}$, we can see that \eq{3OC} holds. Using the second equality in \eq{ids}, this is equivalent to \eq{pS}.
\end{proof}

\Lem{4OC}{
	Let $f,g: I\to\R$ be four times differentiable functions on $I$ with nonvanishing first derivatives. Let $p=(p_1,p_2): I\to\R_{+}^2$ and $q=(q_1,q_2): I\to\R_{+}^2$ be twice differentiable functions such that \eq{2OC} and \eq{3OC} hold on $I$. If $\partial_1^2\partial_2^2 A_{f,p}+\partial_1^2\partial_2^2 A_{g,q}=0$ holds on $\diag(I^2)$, then
	\Eq{4OC}{
	 S(f)'-2\bigg(\frac{(p_1p_2)'}{p_1p_2}+\frac{f''}{f'}\bigg)S(f)
	 +S(g)'-2\bigg(\frac{(q_1q_2)'}{q_1q_2}+\frac{g''}{g'}\bigg)S(g)=0
    }
	is valid on $I$. Consequently, 
	\Eq{pS+}{
	  (p_1-p_2)(S(f)'+S(g)')=0
	}
	is valid on $I$.
}

\begin{proof} The equality $(\partial_1^2\partial_2^2 A_{f,p}+\partial_1^2\partial_2^2 A_{g,q})\circ\Delta_2=0$, formula (iv) of \thm{DB} for $\partial_1^2\partial_2^2 A_{f,p}\circ\Delta_2$ and the analogous expression for $\partial_1^2\partial_2^2 A_{g,q}\circ\Delta_2$ imply that
\Eq{*}{
0 =&\bigg(\Big(\frac{p_1p_2}{p_0^2}\Big)''
    +\frac38\Big(6-\frac{p_0^2}{p_1p_2}\Big)\bigg(\Big(\frac{p_1-p_2}{p_0}\Big)'\bigg)^2\bigg)\bigg(\frac{(p_1p_2)'}{p_1p_2}+\frac{f''}{f'}\bigg)\\
   &-\Big(\frac{p_1p_2}{p_0^2}\Big)'\bigg(\frac{(p_1p_2)'}{p_1p_2}+\frac{f''}{f'}\bigg)'
   -\frac12\Big(\frac{p_1p_2}{p_0^2}\Big)^2\bigg(6-\frac{p_0^2}{p_1p_2}\bigg)\bigg(\frac{(p_1p_2)'}{p_1p_2}+\frac{f''}{f'}\bigg)^3\\
   &+\Big(\frac{p_1p_2}{p_0^2}\Big)^2\bigg(6-\frac{p_0^2}{p_1p_2}\bigg)\bigg(\frac{(p_1p_2)'}{p_1p_2}+\frac{f''}{f'}\bigg)S(f)
   -\Big(\frac{p_1p_2}{p_0^2}\Big)^2S(f)'\\
    &+\bigg(\Big(\frac{q_1q_2}{q_0^2}\Big)''
    +\frac38\Big(6-\frac{q_0^2}{q_1q_2}\Big)\bigg(\Big(\frac{q_1-q_2}{q_0}\Big)'\bigg)^2\bigg)\bigg(\frac{(q_1q_2)'}{q_1q_2}+\frac{g''}{g'}\bigg)\\
   &-\Big(\frac{q_1q_2}{q_0^2}\Big)'\bigg(\frac{(q_1q_2)'}{q_1q_2}+\frac{g''}{g'}\bigg)'
   -\frac12\Big(\frac{q_1q_2}{q_0^2}\Big)^2\bigg(6-\frac{q_0^2}{q_1q_2}\bigg)\bigg(\frac{(q_1q_2)'}{q_1q_2}+\frac{g''}{g'}\bigg)^3\\
   &+\Big(\frac{q_1q_2}{q_0^2}\Big)^2\bigg(6-\frac{q_0^2}{q_1q_2}\bigg)\bigg(\frac{(q_1q_2)'}{q_1q_2}+\frac{g''}{g'}\bigg)S(g)
   -\Big(\frac{q_1q_2}{q_0^2}\Big)^2S(g)'.
}
Using the identities in \eq{ids}, this equality reduces to
\Eq{*}{
	0&=\Big(\frac{p_1p_2}{p_0^2}\Big)^2\bigg(S(f)'+\bigg(\frac{p_0^2}{p_1p_2}-6\bigg)\bigg(\frac{(p_1p_2)'}{p_1p_2}+\frac{f''}{f'}\bigg)S(f)\bigg)\\
   &\qquad+\Big(\frac{q_1q_2}{q_0^2}\Big)^2\bigg(S(g)'+\bigg(\frac{q_0^2}{q_1q_2}-6\bigg)\bigg(\frac{(q_1q_2)'}{q_1q_2}+\frac{g''}{g'}\bigg)S(g)\bigg),
}
which, by the first equality in \eq{ids} yields
\Eq{abc}{
  0&=S(f)'+\bigg(\frac{(p_1-p_2)^2}{p_1p_2}-2\bigg)\bigg(\frac{(p_1p_2)'}{p_1p_2}+\frac{f''}{f'}\bigg)S(f)\\
  &\qquad+S(g)'+\bigg(\frac{(q_1-q_2)^2}{q_1q_2}-2\bigg)\bigg(\frac{(q_1q_2)'}{q_1q_2}+\frac{g''}{g'}\bigg)S(g).
}
On the other hand, by \eq{3OC}, we have
\Eq{*}{
	\frac{p_1-p_2}{p_0}S(f)&=\frac{q_2-q_1}{q_0}S(g),
}
whence we can obtain
\Eq{*}{
  &\frac{(p_1-p_2)^2}{p_1p_2}\bigg(\frac{(p_1p_2)'}{p_1p_2}+\frac{f''}{f'}\bigg)S(f)
  =\frac{p_0^2}{p_1p_2}\cdot\frac{p_1-p_2}{p_0}\bigg(\frac{(p_1p_2)'}{p_1p_2}+\frac{f''}{f'}\bigg)\bigg(\frac{p_1-p_2}{p_0}S(f)\bigg)\\
  &\qquad=-\frac{q_0^2}{q_1q_2}\cdot\frac{q_2-q_1}{q_0}\bigg(\frac{(q_1q_2)'}{q_1q_2}+\frac{g''}{g'}\bigg)\bigg(\frac{q_2-q_1}{q_0}S(g)\bigg)
  =-\frac{(q_1-q_2)^2}{q_1q_2}\bigg(\frac{(q_1q_2)'}{q_1q_2}+\frac{g''}{g'}\bigg)S(g).
}
Using this equality, equation \eq{abc} reduces to \eq{4OC}.

Finally, multiplying \eq{4OC} by $(p_1-p_2)$ side by side and use \eq{2OC} and then \eq{pS}, the last assertion, i.e., equality \eq{pS+} follows.
\end{proof}

\Thm{Main}{
	Let $f,g: I\to\R$ be four times continuously differentiable functions on $I$ with nonvanishing first derivatives. Let $p=(p_1,p_2): I\to\R_{+}^2$ be a twice continuously differentiable function and $q=(q_1,q_2): I\to\R_{+}^2$. Assume that the set 
	\Eq{*}{
	  P:=\{x\in I\mid p_1(x)= p_2(x)\}
	}
	is nowhere dense in $I$. Then the following assertions are equivalent to each other.
	\begin{enumerate}[(i)]
	 \item The invariance equation \eq{inv} holds for every $(x,y)\in I^2$.
	 \item There exists an open set $U\subseteq I^2$ containing the diagonal $\diag(I^2)$ such that the invariance equation \eq{inv} holds for all $(x,y)\in U$.
	 \item The function $q=(q_1,q_2): I\to\R_{+}^2$ is twice continuously differentiable and the system of equalities 
	 \Eq{sys}{
        \partial_1 A_{f,p}+\partial_1 A_{g,q}&=1,\\
        \partial_1\partial_2 A_{f,p}+\partial_1\partial_2 A_{g,q}&=0,\\
        \partial_1^2\partial_2 A_{f,p}+\partial_1^2\partial_2 A_{g,q}&=0,\\ \partial_1^2\partial_2^2 A_{f,p}+\partial_1^2\partial_2^2 A_{g,q}&=0
    }
    holds on the diagonal $\diag(I^2)$.
    \item There exists a real constant $\gamma\in\R$, there exist solutions $u,v,w,z:I\to\R$ of the second-order linear differential equation $F''=\gamma F$ such that $v>0$ and $z>0$ holds on $I$ and $\{u,v\}$ and $\{w,z\}$ are linearly independent such that \eq{fgpq} holds.
	\end{enumerate}
}

\begin{proof} The implication (i)$\Rightarrow$(ii) is trivial. 

Now assume that (ii) is valid. Rearranging the invariance equation \eq{inv}, we get
\Eq{*}{
  A_{g,q}(x,y)=x+y-A_{f,p}(x,y)=:M(x,y) \qquad ((x,y)\in U).
}
By the regularity assumptions on $f$ and $p$, it follows that the mean $M$ defined by the above equality is twice continuously partially differentiable on $U$. Therefore, \lem{reg} implies that $q_1$ and $q_2$ are also twice continuously differentiable on $I$. In view of \thm{DB}, we now obtain that the partial derivatives $\partial_1^i\partial_2^j A_{f,p}$ and $\partial_1^i\partial_2^jA_{g,q}$ exist on $\diag(I^2)$ for all $i,j\in\{1,2\}$. Differentiating both sides of the invariance equation \eq{inv} partially, $i$ and $j$ times with respect to the variables $x$ and $y$, we obtain that the system of equalities \eq{sys} holds on $\diag(I^2)$. This proves that (iii) follows from (ii).

Assume that (iii) is valid. Using \lem{3OC} and the third equality in \eq{sys}, it follows that \eq{pS} is valid on $I$. Therefore, for all $x\in I\setminus P$, we have that
\Eq{Sf=Sg}{
  S(f)(x)=S(g)(x).
}
Observe that $S(f)$ and $S(g)$ are continuous functions, hence this equality is also valid on the closure of $I\setminus P$, which equals $I$ since the set $P$ is nowhere dense. In view of \lem{3OC}, the fourth equality in \eq{sys} implies that \eq{pS+} is valid on $I$. Therefore, by the nowhere density of the set $P$ again, 
\Eq{*}{
  S(f)'(x)+S(g)'(x)=0.
}
Differentiating \eq{Sf=Sg} side by side, the last equation implies that 
\Eq{*}{
  S(f)'(x)=S(g)'(x)=0.
}
Therefore $S(f)$ and $S(g)$ are constant functions which are equal to each other. Let us denote this constant value by $-2\gamma$. Then, applying \lem{uv}, we can conclude that there exist solutions $u,v,w,z:I\to\R$ of the second-order linear differential equation $F''=\gamma F$ such that $v>0$ and $z>0$ holds on $I$ and the first two equalities in \eq{fgpq} are satisfied, i.e., $f=u/v$ and $g=w/z$. The strict monotonicity of $f$ and $g$ imply that $\{u,v\}$ and $\{w,z\}$ are also linearly independent. Observe that
\Eq{*}{
  (u'v-uv')'=u''v+u'v'-u'v'-uv''=\gamma uv-\gamma uv=0.
}
This implies that $u'v-uv'=\alpha$ and, analogously, $w'z-wz'=\beta$ for some nonzero real constants $\alpha$ and $\beta$. Thus
\Eq{*}{
  f'=\frac{u'v-uv'}{v^2}=\frac{\alpha}{v^2} \qquad\mbox{and}\qquad
  g'=\frac{w'z-wz'}{z^2}=\frac{\beta}{z^2}.
}
On the other hand, by \lem{1OC} and \lem{2OC}, the first two equalities in \eq{sys} imply that there exists a nonzero constant $\delta$ such that \eq{ga} is valid on $I$. Consequently,
\Eq{*}{
 p_1q_1=p_2q_2
 =\sqrt{\frac{\delta}{f'g'}}
 =\sqrt{\frac{\delta v^2 z^2}{\alpha\beta}}
 =\sqrt{\frac{\delta }{\alpha\beta}}vz=\eta vz,
}
where $\eta:=\sqrt{\frac{\delta}{\alpha\beta}}$ is a nonzero constant. Define
$\bar{u}:=\eta u$ and $\bar{v}:=\eta v$. Then we can see that $\bar{u}$ and $\bar{v} $ are also solutions of the second-order linear differential equation $F''=\gamma F$ and \eq{fgpq} is satisfied if we replace $u$ and $v$ by $\bar{u}$ and $\bar{v}$, respectively. Hence we have proved that assertion (iii) implies statement (iv).

The implication (iv)$\Rightarrow$(i) is a consequence of \thm{Suff}.
\end{proof}

\def\MR#1{}


\providecommand{\bysame}{\leavevmode\hbox to3em{\hrulefill}\thinspace}
\providecommand{\MR}{\relax\ifhmode\unskip\space\fi MR }
\providecommand{\MRhref}[2]{%
  \href{http://www.ams.org/mathscinet-getitem?mr=#1}{#2}
}
\providecommand{\href}[2]{#2}

\end{document}